
\input amssym.def   
\input amssym.tex   


\magnification=\magstep1

\def \titlefont{\bf}
\font \authorfont=cmcsc10
\def \thanks#1{\footnote{}{{\rm #1}}}

\def \qedbox{{\vbox{\hrule\hbox{\vrule\phantom{q}\vrule}\hrule}}}
\def \qed{\qquad\qedbox\hfill\bigbreak}
\def \cite#1{{\bf [#1]}}
\def \say#1.
{\medbreak\noindent{\it #1.}\enspace}
\def \endsay{\medbreak}


\def \alp{\alpha}
\def \gam{\gamma}
\def \lam{\lambda}
\def \Lam{\Lambda}
\def \QQ{{\Bbb Q}}
\def \ZZ{{\Bbb Z}}
\def \cU{{\cal U}}
\def \cL{{\cal L}}
\def \cO{{\cal O}}
\def \cR{{\cal R}}
\def \zehat{\hat\zeta}
\def \sub{\subset}
\def \tsr{\otimes}
\def \dual{^\vee}
\def \isom{\cong}
\def \End{\mathop{\rm End}\nolimits}
\def \into{\hookrightarrow}

\def \MS{{\cal N}}
\def \Jac{J}
\def \HI{H_I}
\def \PI{P_I}

\def \binomial(#1,#2){{\textstyle\Bigl( {{#1}\atop{#2}} \Bigr)}}
\def \oddbin(#1,#2){\left\{ {{#1}\atop{#2}} \right\}}


\centerline{\titlefont ON THE COHOMOLOGY RING OF THE MODULI SPACE}
\centerline{\titlefont OF RANK 2 VECTOR BUNDLES ON A CURVE}
\bigskip
\centerline{\authorfont
A.D. King\thanks{The first author is supported by
 the SERC/EPSRC (grant GR/J38932).}
\& P.E. Newstead\thanks{Both authors are members of
 the VBAC group of Europroj.}}
\bigskip
\centerline{16 February 1995}
\bigskip


\say Abstract.
Let $\MS_g$ be the moduli space of stable holomorphic vector bundles
 of rank 2 and fixed determinant of odd degree over a smooth complex
 projective curve of genus $g$.
This paper proves various properties of the rational cohomology ring
 $H^*(\MS_g)$.
It is shown that the first relation in genus $g$ between the standard
 generators satisfies a recurrence relation in $g$ and that the invariant
 subring for the mapping class group is a complete intersection ring.
(These two results have been obtained independently by Zagier,
 Baranovsky and Siebert \& Tian.)
A Gr\"obner basis is found for the ideal of invariant relations.
A structural formula for $H^*(\MS_g)$ (originally conjectured by
 Mumford) is verified and a natural monomial basis is given.

\say 1991 Mathematics Subject Classification.
\ Primary: 14H60, 14D20.\ Secondary: 14F25.


\beginsection 1. Introduction

Let $C$ be a smooth complex projective curve of genus $g\geq 2$
 and $L_0$ a line bundle of odd degree on $C$.
Let $\MS_C$ be the moduli space of stable holomorphic vector bundles
 of rank 2  and determinant $L_0$ over $C$.
It is a smooth projective variety of complex dimension $3g-3$,
 independent (up to non-canonical isomorphism) of the choice
 of $L_0$.
The celebrated theorem of Narasimhan \& Seshadri \cite{NS} tells us
 that the underlying real manifold depends only on the fundamental
 group of $C$, i.e. on $g$, and we denote this real manifold by $\MS_g$.

This paper is concerned with the structure of the
 rational cohomology ring $H^*(\MS_g)$,
 a subject which was first studied some 30 years ago and which
 has received a considerable amount of recent interest.
(Note: all cohomology in this paper
 has $\QQ$ coefficients, unless otherwise stated.)
As proved by the second author \cite{Ne2},
 this ring is generated by certain `characteristic classes'
 $\alp,\beta,\psi_1,\ldots,\psi_{2g}$.
It carries a natural action of the mapping class group of $C$, which
 factors through the action of $Sp(2g;\ZZ)$---or more naturally
 $Sp(H^1(C;\ZZ))$---which fixes $\alp$ and $\beta$, and for which
 $\{\psi_1,\ldots,\psi_{2g}\}$ is a basis for the fundamental representation.
The invariant subring $\HI^*(\MS_g)$ is generated by $\alp$, $\beta$
 and $\gam$, the invariant quadratic combination of the $\psi_i$'s.
(For the purposes of this paper, one may actually define $\HI^*(\MS_g)$
 to be the subring generated by $\alp$, $\beta$ and $\gamma$,
 with $\gamma$ defined by the formula in \S2.)

Mumford conjectured a complete system of relations,
 which was confirmed by Kirwan \cite{Ki}.
A different method of determining relations was given by Thaddeus \cite{Th}
 who derived an explicit formula for the intersection pairings.
Thaddeus's work made it clear that the structure of the whole ring
 $H^*(\MS_g)$ depends essentially only on that of the invariant subrings
 $\HI^*(\MS_{g'})$ for $g'\leq g$.
In fact, Mumford had also given a conjectural formulation of this
 property (see (2.7) below) at the time he conjectured the relations,
 but this had not passed into general circulation.
We will verify this conjecture (in Proposition 2.5) using
 an argument closely based on \cite{Th}.

Our main result (Theorem 3.1) is that $\HI^*(\MS_g)$ is a complete
 intersection ring determined by three successive relations
 $\zeta_g$, $\zeta_{g+1}$, $\zeta_{g+2}$
 occurring in the sequence generated recursively by
$$
 \zeta_0 = 1
 \qquad\qquad
 \zeta_{n+1} = \alp\zeta_n + n^2 \beta\zeta_{n-1} + 2n(n-1) \gam\zeta_{n-2}.
$$
This result was first discovered in 1991 by Zagier \cite{Za} and subsequently
 independently rediscovered by Baranovsky \cite{Ba}, Siebert \& Tian
 \cite{ST} and the current authors.
The proof we give here differs from the others mentioned, particularly
 in the crucial step of showing that $\zeta_g$ is a relation.
We use the `original' method that was used to determine relations
 in low genus (\cite{Ra}, \cite{Ne2} \S5) and also to give
 partial information in higher genus (\cite{Ne3}).
The surprising new observation is that the partial information so obtained
 is in fact sufficient to show that the first relation in genus
 $g$ is precisely $\zeta_g$.
One may subsequently determine a complete set of relations
 for the invariant subring, by exploiting the fact
 that $\alp$, $\beta$ and $\gam$ are `universal',
 i.e. under any of the $g$
 natural topological embeddings $\MS_{g-1}\into\MS_g$ they
 pull back to the same classes on the lower genus moduli space.
Thus, in particular, any relation between them in genus $g$ is also a
 relation in lower genus.

In the process of proving the main result
 we will (see Proposition 3.3) identify a Gr\"obner basis for the
 ideal of relations satisfied by $\alp$, $\beta$, $\gamma$
 and hence a monomial basis for $\HI^*(\MS_g)$.
This then enables us to extend a result from \cite{Ne1} and \cite{Ne2}
 (recalled in Proposition 2.1) to show that $H^*(\MS_g)$ has the
 monomial basis
$$
 \alpha^i \beta^j \psi_{i_1} \ldots \psi_{i_k}
 \qquad i+k<g,\quad j+k<g,\quad i_1 < \cdots < i_k.
$$
This provides a natural explanation for the multiplicities
 of the $Sp(2g,\ZZ)$-representations $\Lam^kH^3$ in $H^*(\MS_g)$,
 which were computed by Nelson \cite{N}.

The paper is laid out as follows.
In Section 2, we recall various known and essentially known results,
 presented in the form needed later.
In Section 3, we describe the main new results, in particular the recursively
 generated sequence $(\zeta_g)$ and the Gr\"obner basis.
Sections 4 \& 5 contain the proofs.


\say Acknowledgement.
This work arose out of joint work with V. Balaji.
We are grateful to him for helping to set the ball rolling
 and to C.T.C. Wall for some remarks that kept it on course.
We would also like to thank V. Baranovsky, B. Siebert \& G. Tian,
 and D. Zagier for discussing their results with us.
Finally, we should mention that, even though they do not play a r\^ole
 in the final proof, Thaddeus' formula for the intersection pairings
 and MAPLE were invaluable help in finding the results described here.


\beginsection 2. Preliminary results.

In this section, we bring together various results
 which are relevant to this paper.
In some cases, the formulation of these results may not
 be so familiar.

The Betti numbers of $\MS_g$ were calculated in \cite{Ne1} and
 may be efficiently presented as the coefficients of the
 Poincar\'e polynomial (\cite{Ha}, Satz 3.3)
$$
 P(\MS_g;t) = {(1+t^3)^{2g} - t^{2g}(1+t)^{2g}\over (1-t^2)(1-t^4)}.
 \eqno{(2.1)}
$$
This may easily be interpreted as the difference of two expressions,
 the first of which enumerates the free ring on a set of generators of
 $H^*(\MS_g)$ and the second of which enumerates the ideal of relations.

A suitable set of generators was calculated in \cite{Ne2}.
It consists of classes $\alp\in H^2$, $\beta\in H^4$ and
 a basis $\psi_1,\ldots,\psi_{2g}$ for $H^3$.
These may be defined by splitting the second Chern class of
 the universal endomorphism bundle on $\MS_C\times C$ into its
 K\"unneth components as follows
$$
 c_2(\End \cU) = 2\alp f + 4\psi - \beta
 \qquad\qquad
 \psi=\sum_{i=1}^{2g}\psi_i e_i,
 \eqno{(2.2)}
$$
 where $f$ is the positive generator of $H^2(C;\ZZ)$,
 and $e_1,\ldots,e_{2g}$ is a symplectic basis of $H^1(C;\ZZ)$
 with $e_ie_{i+g} =-f$ for $1\leq i \leq g$.
Note that, because odd degree classes in cohomology anti-commute,
 the `free ring on the generators' means the ring
 $\Lam^*H^3\tsr\QQ[\alp,\beta]$, whose Poincar\'e polynomial is
$$
 {(1+t^3)^{2g}\over (1-t^2)(1-t^4)}.
$$
The additional basic invariant class $\gam$ is given by
$$
 \gam = \int_C \psi^2 = 2\sum_{i=1}^g \psi_i \psi_{i+g}.
$$

Now, results from \cite{Ne1} and \cite{Ne2} suffice to identify
 a monomial basis for $H^*(\MS_g)$.

\proclaim Proposition 2.1.
If $n\leq 3g-1$, then the monomials
$$
 \alp^i\beta^j\psi_{i_1}\ldots\psi_{i_k}
 \qquad 2i+4j+3k=n,
 \quad i+k < g,\quad i_1 < \cdots < i_k
 \eqno{(2.3)}
$$
 form a basis of $H^n(\MS_g)$.

\say Proof.
The fact that these monomials are independent follows from
 Propositions 2.6 \& 3.4 of \cite{Ne2},
 together with Lemma 4 of \cite{Ne1} which says that multiplication
 by $\beta$ is injective on $H^{n-4}(\MS_g)$ for $n\leq 3g-1$.
Since the number of such monomials is equal to the $n$th Betti number
 (\cite{Ne1} Theorem 2) the result follows.
\qed

\say Remark 2.2.
One may extend the monomial basis (2.3) to one for the whole
 cohomology ring by exploiting the fact that $\alp$ is ample,
 and hence the Hard Lefschetz Theorem implies that, for $n\leq 3g-3$,
 multiplication by $\alp^{3g-3-n}$ provides an isomorphism
 $H^n(\MS_g)\isom H^{6g-6-n}(\MS_g)$.
However, in Proposition 3.3(iii) below, we will find a more
 natural extension.

\say Remark 2.3.
Proposition 2.1 leads one to express the Poincar\'e polynomial
 of $\MS_g$ as
$$
 \sum_{k=0}^{g-1} \binomial(2g,k) t^{3k}
 {(1-t^{2g-2k})(1-t^{4g-4k})\over(1-t^2)(1-t^4)},
 \eqno{(2.4)}
$$
 because the coefficient of $t^n$ counts the monomials in (2.3)
 for $n\leq 3g-3$ (when the term $(1-t^{4g-4k})$ is irrelevant),
 while the whole expression is
 compatible with Poincar\'e duality.
Now, (2.4) can be equated with (2.1), or equivalently
$$
 \sum_{k=0}^{2g}  \binomial(2g,k)
 {t^{3k} - t^{2g+k} \over (1-t^2)(1-t^4)},
 \eqno{(2.5)}
$$
 using only the fact that $\binomial(2g,k)=\binomial(2g,2g-k)$.
This observation is important because it also applies
 to the calculation of the Poincar\'e polynomial for the algebraic
 cohomology in \cite{BKN}.
\endsay

A suitable set of relations in $\Lam^*H^3\tsr\QQ[\alp,\beta]$
 was proposed by Mumford and the
 fact that it is a complete set was proved by Kirwan \cite{Ki}.
Briefly, these relations are obtained from a natural rank $2g-1$
 vector bundle $F$ over $\MS_C\times\Jac$, where $\Jac$ is the
 Jacobian of $C$.
The Chern classes of $F$ may be expressed in terms of
 $\alp,\beta,\psi_1,\ldots,\psi_{2g}$ and the cohomology classes of $\Jac$.
Hence, for each class $\omega\in H^*(\Jac)$,
 one may formally evaluate
$$
 \int_{\Jac} c_r(F)\omega
$$
 to obtain polynomials in $\alp,\beta,\psi_1,\ldots,\psi_{2g}$.
When $r\geq 2g$, these must be relations, because $c_r(F)=0$.
These are the Mumford relations.

\say Remark 2.4.
Observe that the relation ideal must have Poincar\'e polynomial
$$
 {t^{2g}(1+t)^{2g}\over (1-t^2)(1-t^4)}
$$
 and that $(1+t)^{2g}=P(\Jac;t)$.
These facts strongly suggest that the ideal
 is freely generated as a $\QQ[\alp,\beta]$-module
 by the Mumford relations coming just from $c_{2g}(F)$.
\endsay

Mumford actually made an additional (but less well publicised)
 conjecture about the relations,
 based on Bayer's computer calculations for $g\leq5$.
For $g\geq2$ let
$$
 I_g\sub\QQ[\alp,\beta,\gam]
$$
 be the ideal of relations holding in $H_I^*(\MS_g)$,
 and in addition let $I_0$=(1) and $I_1=(\alp,\beta,\gam)$.
The `universality' of $\alp,\beta,\gam$ (c.f. \S1) implies
 that $I_{g+1}\sub I_g$ for all $g$.
For $0\leq k\leq g$, recall that the primitive component of
 $\Lam^kH^3$ is defined by
$$
 \Lam_0^k H^3 =
 \ker \bigl( \gam^{g-k+1}:\Lam^k H^3 \to \Lam^{2g-k+2} H^3 \bigr).
$$
(Note: $\gam^{g-k+1}$ is acting by multiplication in the
 exterior algebra $\Lam^*H^3$, not in $H^*(\MS_g)$.)
This gives rise to the primitive decomposition of the exterior algebra
$$
 \Lam^*H^3 \isom \bigoplus_{k=0}^g \Lam_0^kH^3
 \tsr \QQ[\gam]/\bigl(\gam^{g-k+1}\bigr).
 \eqno{(2.6)}
$$
Mumford conjectured that a similar decomposition holds for
 the cohomology ring, in that there is an isomorphism
$$
 H^*(\MS_g) \isom \bigoplus_{k=0}^g
 \Lam_0^kH^3(\MS_g) \otimes \QQ[\alp,\beta,\gam]/I_{g-k}.
 \eqno{(2.7)}
$$
Note that the general form of the expression is to be expected from
 the way the Mumford relations are generated;
it is the fact that the annihilator of $\Lam_0^k H^3(\MS_g)$ depends
 only on $g-k$ that gives this second conjecture its extra precision.
The basic properties of the intersection pairings, proved in \cite{Th},
 can be used to verify the conjecture.

\proclaim Proposition 2.5.
Let
$$
 \nu: \bigoplus_{k=0}^g
 \Lam_0^kH^3(\MS_g) \otimes \QQ[\alp,\beta,\gam]
 \to H^*(\MS_g)
$$
 be the obvious map.
Then
$$
 \ker\nu = \bigoplus_{k=0}^g \Lam_0^kH^3(\MS_g) \otimes I_{g-k}
$$

\say Proof.
There is a coarser decomposition than
 the primitive decomposition, namely that when $k$ is even
$$
 \Lam^kH^3 = \ker\bigl(\gam^{g-{k\over2}}\bigr)
 \oplus \bigl\langle \gam^{{k\over2}} \bigr\rangle.
$$
One may render \cite{Th} Proposition 24, et seq. as follows:
 if $p\in\QQ[\alp,\beta,\gam]$ and $q\in\Lam^k H^3$,
 and either $k$ is odd or $q\in \ker\bigl(\gam^{g-{k\over2}}\bigr)$,
 then
$$
 \int_{\MS_g} pq = 0.
$$
Now observe that if $q\in\Lam^k_0 H^3$
 and $q'\in\Lam^{k'}_0 H^3$, with $k+k'$ even
 and $k\neq k'$, then
 $qq'\in \ker\bigl(\gam^{g-{k+k'\over2}}\bigr)$.
Further observe that the pairing
$$
 \Lam^k_0 H^3 \tsr \Lam^k_0 H^3 \to \Lam^{2g} H^3:
 (q,q') \mapsto qq'\gam^{g-k}
$$
 is non-degenerate.
Hence we may, for $0\leq k\leq g$, choose a basis
 $\bigl\{ q_{ki} \mid 1\leq i \leq \dim \Lam^k_0H^3 \bigr\}$
 for $\Lam_0^k H^3$,
 and also a dual basis $\bigl\{ q_{ki}\dual \bigr\}$
 with respect to this pairing.

Consider a general element
$$
 r = \sum_{l,j} q_{lj} p_{lj}
 \in  \bigoplus_{k=0}^g \Lam_0^kH^3 \otimes \QQ[\alp,\beta,\gam]
$$
 where $p_{lj}\in \QQ[\alp,\beta,\gam]$.
Then the above observations, together with
 \cite{Th} Proposition 26, show that, for $p'\in\QQ[\alp,\beta,\gam]$,
$$
 \int_{\MS_g} r q_{ki}\dual p'
 \quad = \quad
 \int_{\MS_g} p_{ki}p' q_{ki} q_{ki}\dual
 \quad\propto\quad
 \int_{\MS_g} p_{ki}p' \gam^k
 \quad\propto\quad
 \int_{\MS_{g-k}} p_{ki}p'.
$$
This evaluation is zero for all $p'$ if and only if
 $p_{ki}\in I_{g-k}$.
Hence $r\in\ker\nu$ if and only if $p_{ki}\in I_{g-k}$
 for all $k$ and $i$.
\qed

{}From this result we may calculate the Poincar\'e polynomial
 $\PI(\MS_g;t)$ of $\HI^*(\MS_g)$.
{}From (2.7) and (2.4) we have
$$
 \sum_{k=0}^{g-1}
 \left( \binomial(2g,k)-\binomial(2g,k-2) \right)
 t^{3k} \PI(\MS_{g-k};t)
 =
 \sum_{k=0}^{g-1} \binomial(2g,k) t^{3k}
{(1-t^{2g-2k})(1-t^{4g-4k})\over(1-t^2)(1-t^4)}
$$
 from which we obtain by induction
$$
 \PI(\MS_g;t) - t^6 \PI(\MS_{g-2};t)
 =
{(1-t^{2g})(1-t^{4g})\over(1-t^2)(1-t^4)}
$$
 and hence
$$
 \eqalign{
 \PI(\MS_g;t)
 &=
 \sum_{p=0}^{[{g\over 2}]}
 {\bigl(1-t^{2g-4p}\bigr)\bigl(1-t^{4g-8p}\bigr)
 \over (1-t^2)(1-t^4)} t^{6p} \cr
 &=
 { \bigl( 1-t^{2g} \bigr) \bigl( 1-t^{2g+2} \bigr) \bigl( 1-t^{2g+4} \bigr)
 \over \bigl( 1-t^2 \bigr) \bigl( 1-t^4 \bigr) \bigl( 1-t^6 \bigr) } }
 \eqno{(2.8)}
$$

This formula was the starting point for this paper,
 when it appeared in the course of the proof that $H_I^*(\MS_g)$
 is the algebraic cohomology ring for a general curve $C$ (see \cite{BKN}).


\beginsection 3. Main results.

The fact that (2.8) is the Poincar\'e polynomial of $H_I^*(\MS_g)$
 strongly suggests that this ring is a complete intersection,
 i.e. that $I_g$ is generated by three relations
 in degrees $2g$, $2g+2$ and $2g+4$.
In this section, we present the results which show that this is
 indeed the case and also describe a Gr\"obner basis for $I_g$.
The proofs appear in the following two sections.

Consider the sequence of homogeneous polynomials in $\ZZ[\alp,\beta,\gam]$
 defined recursively by
$$
\zeta_{n+1} = \alp \zeta_n + n^2 \beta \zeta_{n-1} + 2n(n-1) \gam \zeta_{n-2}
\eqno{(3.1)}
$$
starting with $\zeta_0=1$.
Observe that one conveniently does not need to specify the values of
 $\zeta_{-1}$ and $\zeta_{-2}$.

One may partially solve this recurrence formula to write
$$
 \zeta_g = \sum_{{s,t,u\geq0 \atop s+2t+3u=g}}
 \lam_{t,u} \binomial(g,s) \alp^s \beta^t (2\gam+\alp\beta)^u
 \eqno{(3.2)}
$$
 where the constants $\lam_{t,u}$ are independent of $g$ and
 are in turn defined recursively by
$$
 \lam_{t,u} = (d-1)^2 \lam_{t-1,u} + (d-1)(d-2) \lam_{t,u-1},
 \eqno{(3.3)}
$$
 where $d=2t+3u$, starting with $\lam_{0,0}=1$
 and $\lam_{t,u}=0$ if $t<0$ or $u<0$.
The first few values are
$$
 \lam_{1,0}=1 \quad \lam_{0,1}=2 \quad \lam_{2,0}=9 \quad
 \lam_{1,1}=44 \quad \lam_{0,2}=40 \quad \lam_{3,0}=225
$$
One may also observe that (3.2) agrees with the formula for the
relations modulo $\beta$ given in \cite{Ne3}.

The main theorem, proved in \cite{Za}, \cite{Ba}, \cite{ST} and here,
 is

\proclaim Theorem 3.1.
Let $I_g\sub \QQ[\alp,\beta,\gam]$ and $\zeta_g$ be as defined above
 ($I_g$ in \S2).
\item{i)} $\zeta_g\in I_g$ for all $g$,
\item{ii)} $I_g = (\zeta_g, \zeta_{g+1}, \zeta_{g+2})$.

Of these, the first part requires the main geometric input and has
 the widest variety of proofs.
In \cite{Za}, it is proved by showing that
 $\zeta_n$ is a Mumford relation for $n\geq g$,
 while in \cite{Ba} and \cite{ST} it follows from the fact that
 $\zeta_n/n!$ is the $n$th Chern class of a rank $g-1$ bundle.
In this paper we take a rather different approach.
We will show that the `first relation in genus $g$' (known to be
 in degree $2g$) agrees with $\zeta_g$ for small $g$ and also
 satisfies the recurrence formula (3.1),
 hence it is $\zeta_g$ for all $g$.

Note that the first part immediately implies that
 $\zeta_n\in I_g$ for all $n\geq g$.
The second part requires some algebraic computation and
 a counting argument.
We do this by finding an explicit Gr\"obner basis for
 the ideal $(\zeta_g,\zeta_{g+1},\zeta_{g+2})$
 and using (2.8).
The determination of this basis rests on the following explicit
 calculation.

\proclaim Lemma 3.2.
Letting $\zehat_n=\zeta_n/n!$, the following definitions are equivalent
$$
 \eqalignno{
 \zehat_{g,n} &= \sum_{i=0}^n
 {1\over i!} \binomial(g-i,n-i) (2\gam)^i\beta^{n-i} \zehat_{g-n-i}
 & (3.4a) \cr
 (-1)^n\zehat_{g,n} &= \sum_{i=0}^\infty
 (-1)^i \oddbin(g-n,i) \zehat_{n-i}\zehat_{g+i}
 & (3.4b) }
$$
 where we use the (non-standard) notation
 $\oddbin(k,i)=\binomial(k+i,i) + \binomial(k+i-1,i-1)$
 and the convention that $\zehat_n=0$ for $n<0$.

When $g$ is odd, we also let
$$
 \zehat_{g,{g\over2}} = 6\gam \zehat_{g,{g-3\over2}}
 - {g-1\over 4}\alp^2 \zehat_{g,{g-1\over2}},
 \eqno{(3.5)}
$$
 which is chosen to have leading monomial $\gam^{g+1\over2}$.
We may then prove the following.

\proclaim Proposition 3.3.
\item{i)} $\zehat_{g,0},\ldots,\zehat_{g,g}$,
 together with $\zehat_{g,{g\over2}}$ when $g$ is odd, is a Gr\"obner basis
 for $I_g$.
\item{ii)} $H_I^*(\MS_g)$ has a monomial basis consisting of
$$
 \alp^i\beta^j\gam^p \qquad i+2p<g,\quad j+2p<g.
 \eqno{(3.6)}
$$
\item{iii)} $H^*(\MS_g)$ has a monomial basis consisting of
$$
 \alpha^i \beta^j \psi_{i_1} \ldots \psi_{i_k}
 \qquad i+k<g,\quad j+k<g,\quad i_1 < \cdots < i_k.
 \eqno{(3.7)}
$$

Notice that $\zehat_{g,g}=\beta^g$.
The fact that this is a relation was conjectured in
 \cite{Ne2} and proved in \cite{Ki} and \cite{Th}.
Note also that the monomials in (3.7) are naturally counted
 by (2.4).


\beginsection 4. Proof of Theorem 3.1(i).

Note first that the cases $g=0,1$ follow by direct computation,
 so we may suppose that $g\geq 2$.

By Proposition 2.1 (see also \cite{Ne2} \S5), the first relation in
 $\Lam^*H^3(\MS_g)\tsr \QQ[\alp,\beta]$ is in degree $2g$
 and is unique up to scalar multiplication.
In addition, the coefficient of $\alp^g$ is non-zero,
 so we normalise the relation by choosing this coefficient to be 1,
 and we denote this normalised relation by $r_g$.
By Proposition 2.5, the uniqueness implies that
 $r_g\in\QQ[\alp,\beta,\gamma]$ and so we may write it
$$
 r_g = \sum_{{s,t,u\geq0\atop s+2t+3u=g}}
 \mu_{t,u} \alp^s\beta^t(2\gam+\alp\beta)^u
$$
 where $\mu_{t,u}$ depends on $g$ as well as $(t,u)$.

The `original' method for obtaining partial information about
 the relations was to restrict to the base of a special family of bundles
 whose cohomology was known.
This was how Ramanan \cite{Ra} found the complete
 set of relations in genus 3, and we shall use the same family here.
More precisely, we obtain the following information about $r_g$,
 which turns out to be sufficient to prove that $r_g=\zeta_g$.

\proclaim Proposition 4.1.
For all $d\geq 0$,
$$
 \sum_{{t,u\geq0\atop 2t+3u=d}} \mu_{t,u}
 = \binomial(g,d) \sum_{j=0}^d (-1)^j {d!\over j!}
$$

We shall prove this formula at the end of the section.
In fact, we will need only the first 5 non-trivial cases, namely
$$
 \eqalign{
 \mu_{1,0} = \binomial(g,2) \quad
 \mu_{0,1} = 2\binomial(g,3)& \quad
 \mu_{2,0} = 9\binomial(g,4) \quad
 \mu_{1,1} = 44\binomial(g,5) \cr
 \mu_{0,2} + \mu_{3,0} &= 265\binomial(g,6) }
$$
To see that this is sufficient, first observe that for $g\leq 5$
 all the coefficients are determined and agree with those of $\zeta_g$
 as given by (3.2).
Hence $r_g=\zeta_g$ for $g\leq 5$.

For $g\geq 5$, suppose we know that $r_m=\zeta_m$ for $m\leq g$,
 and hence that $\zeta_m\in I_n$ for $n\leq m\leq g$.
Now $r_{g+1}\in I_{g+1}\sub I_{g-2}$.
However, Corollary 2.6 shows that there are
 at most 6 independent relations in $I_{g-2}$ of degree $2g+2$
 and we already know the 6 relations
$$
 \alp \zeta_g,\quad
 \alp^2 \zeta_{g-1},\quad
 \alp^3 \zeta_{g-2},\quad
 (2\gam+\alp\beta)\zeta_{g-2},\quad
 \beta \zeta_{g-1},\quad
 \alp\beta \zeta_{g-2}.
 \eqno{(4.1)}
$$
The following matrix calculates the (0,0), (0,1), (1,0), (2,0)
 and (1,1) coefficients and the sum of the (0,2) and (3,0)
 coefficients of a linear combination of the above 6 relations.
(Note: the ``$(t,u)$ coefficient'' means
 the coefficient of $\alp^s\beta^t(2\gam+\alp\beta)^u$,
 with $s+2t+3u=g+1$.)
\medskip
$$
 \pmatrix{\cr
 1 & 1 & 1 & 0 & 0 & 0 \cr\cr
 2 \binomial(g, 3) & 2 \binomial(g - 1, 3) & 2 \binomial(g - 2, 3) &
 1 & 0 & 0 \cr\cr
 \binomial(g, 2) & \binomial(g - 1, 2) & \binomial(g - 2, 2)
 & 0 & 1 & 1 \cr\cr
 9 \binomial(g, 4) & 9 \binomial(g - 1, 4) & 9 \binomial(g - 2, 4) &
 0 & \binomial(g - 1, 2) & \binomial(g - 2, 2) \cr\cr
 44 \binomial(g, 5) & 44 \binomial(g - 1, 5) & 44 \binomial(g - 2, 5) &
 \binomial(g - 2, 2) & 2 \binomial(g - 1, 3) & 2 \binomial(g - 2, 3) \cr\cr
 265 \binomial(g, 6) & 265 \binomial(g - 1, 6) & 265 \binomial(g - 2, 6) &
 2 \binomial(g - 2, 3) & 9 \binomial(g - 1, 4) & 9 \binomial(g - 2, 4) \cr
 &&&&&& }
$$
\medskip\noindent
This matrix has determinant $12 (g - 1) (g - 2)^3 (g - 3)^2 (g - 4)$
 which shows firstly that the 6 relations in (4.1) are independent and
 hence span all the relations in degree $2g+2$, and secondly that
 any relation in this degree is determined by the knowledge of these
 coefficients.
We may use Proposition 4.1 to check that these coefficents
 in $r_{g+1}$ agree with those of $\zeta_{g+1}$ as given by (3.1),
 and hence that $r_{g+1}=\zeta_{g+1}$.

Therefore, to complete the proof of Theorem 3.1(i),
 we need the following.

\say Proof of Proposition 4.1.

The proof goes by constucting an explicit family of bundles parametrised
 by a projective bundle associated to a rank $g$ Picard bundle over
 the Jacobian $\Jac$ and observing that the pull-back of $r_g$
 must be a multiple of the canonical relation (in degree $2g$)
 in the cohomology ring of the projective bundle.
We recall the construction of this family in \cite{Ra}.

For this construction it is convenient to choose $L_0$ to have degree 1 and
 let $\cL$ be the universal bundle on $\Jac\times C$,
 normalised by the requirement that $c_1(\cL)\in H^1(\Jac)\tsr H^1(C)$.
Then, because the degree of $M=\cL^2\tsr L_0^{-1}$ is negative
 along the fibres of $\pi:\Jac\times C \to \Jac$, we see that
 $R^1\pi_*M$ is locally free of rank $g$ (by Riemann-Roch).
Let $P\to\Jac$ be the associated projective bundle and observe
 that there is a universal extension
$$
 0 \to \cL \to E \to \cO_P(-1)\tsr \cL^{-1}\tsr L_0 \to 0
 \eqno{(4.2)}
$$
 over $P\times C$ (where the obvious pull-backs have been omitted
 and $\cO_P(-1)$ is dual to the relative hyperplane bundle).
The restriction of $E$ to $\{p\}\times C$ never splits and hence
 is always a stable bundle of rank 2 and determinant $L_0$.
Hence, there is a morphism $\phi:P\to \MS_C$ such that $E$ is
 the pullback of the universal bundle $\cU$ on $\MS_C\times C$.
{}From (4.2) we may calculate the Chern classes of $E$ to be
 (again omitting pull-backs)
$$
 c_1(E) = f - h \qquad c_2(E) = -h c_1(\cL) - c_1(\cL)^2
$$
 where $h$ is the class of $\cO_P(1)$ and, as before,
 $f$ is the fundamental class of $C$.
In addition, note that $c_1(\cL)^2 = -2\theta f$, where
 $\theta$ is the class of the theta divisor on $\Jac$.
Hence
$$
 c_2(\End E) = 4c_2(E) - c_1(E)^2
 = 2(h+4\theta)f - 4hc_1(\cL) -h^2
$$
 and so, comparing this with (2.2), we see that
$$
 \phi^*(\alp) = h+4\theta \qquad \phi^*(\beta) = h^2
 \qquad \phi^*(2\gam+\alp\beta) = h^3
$$
and thus, writing $\bar\alp$ for $\phi^*(\alp)$, we have
$$
 \phi^*(r_g) = \sum_{{s,t,u\geq0\atop s+2t+3u=g}}
 \mu_{t,u} \bar\alp^s h^{2t+3u}.
$$

On the other hand (c.f \cite{Ra} Lemma 5.4),
 the canonical relation in the cohomology ring of $P$
 (and only relation in degree $2g$ involving just $h$ and $\bar\alp$) is
$$
 \sum_{i=0}^g {(\bar\alp-h)^i\over i!} h^{g-i}
 \eqno{(4.3)}
$$
 because the $i$th Chern class of $R^1\pi_*M$ is $(4\theta)^i/i!$.
Comparing coefficients of $\bar\alp^g$, we see that (4.3) must be
 equal to $\phi^*(r_g)/g!$.
Hence, comparing
 coefficients of $\bar\alp^{g-d} h^d$ we obtain
$$
 \sum_{{t,u\geq0\atop 2t+3u=d}} \mu_{t,u}
 = g! \sum_{j=0}^d {(-1)^j\over (g-d+j)!} \binomial(g-d+j,j)
$$
This is easily rearranged to give the formula of the proposition.
\qed


\beginsection 5. Remaining proofs.

\say Proof of Lemma 3.2.
We prove the equivalence of formulae (3.4a) and (3.4b)
 by showing that both satisfy the recurrence formula
$$
 n\zehat_{g,n} = g\beta \zehat_{g-1,n-1} + 2\gam \zehat_{g-2,n-1}
 \eqno{(5.1)}
$$
 starting with $\zehat_{g,0} = \zehat_g$,
 which is immediate for both formulae.
For (3.4a), the recurrence follows easily, using the identity
$$
 n\binomial(g-i,n-i)
 = g\binomial(g-1-i,n-1-i) + i\binomial(g-1-i,n-i)
$$
For (3.4b), observe first that the sequence $\zehat_m$ satisfies the
 recurrence formula
$$
 (m+1)\zehat_{m+1}
 = \alp \zehat_m + m\beta \zehat_{m-1} + 2\gam \zehat_{m-2}.
$$
Using this to write $\alp\zehat_n\zehat_m $ in two ways, we obtain
$$
 (m+1)\zehat_n\zehat_{m+1} - (n+1)\zehat_{n+1}\zehat_m
 = \beta \bigl( m\zehat_n\zehat_{m-1} - n\zehat_{n-1}\zehat_m \bigr)
 + 2\gam \bigl( \zehat_n\zehat_{m-2} - \zehat_{n-2}\zehat_m \bigr)
$$
In addition, we use the identities
$$
 \eqalign{
 n\oddbin(g-n,i) &= (n-i)\binomial(g-n+i,i) + (g+i)\binomial(g-n+i-1,i-1) \cr
 g\oddbin(g-n,i) &= (g+i)\binomial(g-n+i,i) + (n-i)\binomial(g-n+i-1,i-1) \cr
 \oddbin(g-n-1,i) &= \binomial(g-n+i,i) - \binomial(g-n+i-2,i-2) }
$$
 to write
$$
 \eqalign{
 -n \sum_{i=0}^\infty  &(-1)^i \oddbin(g-n,i) \zehat_{n-i}\zehat_{g+i} \cr
 &= \sum_{i=0}^\infty (-1)^i \binomial(g-n+i,i)
 \bigl( (g+i+1)\zehat_{n-i-1}\zehat_{g+i+1}
 - (n-i)\zehat_{n-i}\zehat_{g+i} \bigr) \cr
 &= \beta \sum_{i=0}^\infty  (-1)^i \binomial(g-n+i,i)
  \bigl( (g+i)\zehat_{n-i-1}\zehat_{g+i-1}
 - (n-i-1)\zehat_{n-i-2}\zehat_{g+i} \bigr) \cr
 &\quad+ 2\gam \sum_{i=0}^\infty (-1)^i \binomial(g-n+i,i)
  \bigl( \zehat_{n-i-1}\zehat_{g+i-2}
 - \zehat_{n-i-3}\zehat_{g+i} \bigr) \cr
 &= g\beta \sum_{i=0}^\infty (-1)^i \oddbin(g-n,i) \zehat_{n-1-i}\zehat_{g-1+i}
 + 2\gam \sum_{i=0}^\infty (-1)^i \oddbin(g-n-1,i) \zehat_{n-1-i}\zehat_{g-2+i}
}$$
which is the required recurrence formula.
\qed

\say Proof of Theorem 3.1(ii).
It is clear from (3.4b) that
 $\zehat_{g,n}\in (\zeta_g,\zeta_{g+1},\zeta_{g+2})$ for all $n$,
 while both equations show that $\zehat_{g,n}=0$ for $n>g$.
{}From (3.4a) one may read off the leading monomials of $\zehat_{g,n}$,
 using the reverse lexographic ordering with the variables ordered
 $\alp$, $\gam$, $\beta$.
They are
$$
 \eqalign{
  \alp^{g-2n}\gam^n &\qquad\hbox{for $0\leq n\leq {g\over 2}$} \cr
  \gam^{g-n}\beta^{2n-g} &\qquad\hbox{for ${g\over2}\leq n\leq g$} }
$$
When $g$ is odd there is no leading monomial which is just a power of
 $\gam$, so we are required to also consider $\zehat_{g,{g\over2}}$,
 as defined in (3.5), which has leading monomial $\gam^{g+1\over2}$.

{}From this it immediately follows that the monomials listed in (3.6) span
 the graded ring
$$
 \cR_g^* = \QQ[\alp,\beta,\gam]/(\zeta_g,\zeta_{g+1},\zeta_{g+2}).
$$
Hence, an upper bound on the Poincar\'e polynomial of $\cR^*_g$ is
$$
 \sum_{p=0}^{[{g\over 2}]}
 {\bigl(1-t^{2g-4p}\bigr)\bigl(1-t^{4g-8p}\bigr)
 \over (1-t^2)(1-t^4)} t^{6p}
$$
But (see (2.8)) this is also the  Poincar\'e polynomial of $H_I^*(\MS_g)$,
 which is a quotient of $\cR_g^*$.
Hence we see that $\cR_g^*=H_I^*(\MS_g)$,
 which proves Theorem 3.1(ii).
\qed

\say Proof of Proposition 3.3.
Part (ii) follows from the fact that the spanning monomials in (3.6) are
 of the correct number and hence a basis, which implies
 part (i), that the set of relations for which these are complementary
 monomials is actually a Gr\"obner basis.
Part (iii) follows from part (ii) and the `structural formula' (2.7).
\qed


\def \paper [#1] #2 (#3) #4; #5 <#6> #7;{
\item{{\bf[#1]}} #2,\ #4,\ {\it #5}\ {\bf #6}\ (#3)\ #7 }

\def \preprint [#1] #2; #3;  #4;{
\item{{\bf[#1]}} #2,\ #3,\ #4 }

\beginsection References.

\preprint [BKN] V. Balaji, A.D. King \& P.E. Newstead;
Algebraic cohomology of the moduli space of
rank 2 vector bundles on a curve;
preprint 1995;

\paper [Ba] V. Baranovsky (1994)
Cohomology ring of the moduli space of stable vector bundles with
odd determinant;
Izv. Russ. Acad. Nauk. <58 n4> 204--210;

\paper [Ha] G. Harder (1970)
Eine Bemerkung zu einer Arbeit von P.E. Newstead;
J. Reine Angew. Math. <242> 16--25;

\paper [Ki] F. Kirwan (1992)
The cohomology ring of moduli spaces of bundles over Riemann surfaces;
J. Amer. Math. Soc. <5> 853--906;

\paper [NS] M.S. Narasimhan \& C.S. Seshadri (1965)
Stable and unitary vector bundles on a compact Riemann surface;
Ann. of Math. <82> 540--567;

\preprint [N] G. Nelson;
The homology of moduli spaces on a Riemann surface
as a representation of the mapping class group;
preprint 1995;

\paper [Ne1] P.E. Newstead (1967)
Topological properties of some spaces of stable bundles;
Topology <6> 241--262;

\paper [Ne2] P.E. Newstead (1972)
Characteristic classes of stable bundles of rank 2 over an
algebraic curve;
Trans. Amer. Math. Soc. <169> 337--345;

\paper [Ne3] P.E. Newstead (1984)
On the relations between characteristic classes of stable bundles
of rank 2 over an algebraic curve;
Bull. Amer. Math. Soc. <10> 292--294;

\paper [Ra] S. Ramanan (1973)
The moduli spaces of vector bundles over an algebraic curve;
Math. Ann. <200> 69--84;

\preprint [ST] B. Siebert \& G. Tian;
Recursive relations for the cohomology ring of
moduli spaces of stable bundles;
preprint alg-geom/9410019;

\paper [Th] M. Thaddeus (1992)
Conformal field theory and the cohomology of the moduli space
of stable bundles;
J. Diff. Geom. <35> 131--149;

\preprint [Za] D. Zagier;
On the cohomology of moduli spaces of rank two vector
bundles over curves;
in preparation since 1991;


\begingroup
\parindent=0pt
\obeylines
\bigskip
Department of Pure Mathematics,
University of Liverpool,
P.O. Box 147,
Liverpool L69 3BX, U.K.
\smallskip
e-mail: aking@liv.ac.uk, newstead@liv.ac.uk
\endgroup

\bye